\documentclass[12pt]{article}
\usepackage{graphicx}
\usepackage{graphics}
\usepackage{amsthm}
\usepackage{amssymb, amsmath}
\title{Reductions of particular hypergeometric functions $_3F_2(a,a+1/3,a+2/3;p/3,q/3;\pm 1)$} %[others later?]
\author{Mark W. Coffey\\
Department of Physics\\
Colorado School of Mines\\
Golden, CO  80401\\
USA\\
mcoffey@mines.edu}
\date{June 6, 2015}
\newcommand*{\Scale}[2][4]{\scalebox{#1}{$#2$}}%
\begin{document}
\maketitle
%\vspace{.25cm}
\baselineskip=25 pt
\begin{abstract}

We principally present reductions of certain generalized hypergeometric functions $_3F_2(\pm 1)$ 
in terms of products of elementary functions.  Most of these results have been known for some time, but one of the methods, wherein we simultaneously solve for three alternating binomial sums, may 
be new.  We obtain a functional equation holding for all three of this set of alternating binomial sums.  Using successive derivatives, we show how related chains of $_3F_2(\pm 1)$ values may be
obtained.  It may be emphasized that we make no reliance on the WZ method for hypergeometric summation.  Additional material on Pochhammer symbols and certain of their products is presented in an Appendix to supplement the pedagogical content of the paper.
% then probably add the $_4F_3(-1)$--check argument--reduction of my Remark of 2012 ...
% also the derivation provided to Victor for another $_3F_2(1)$?
% week of May 25th need to update Abstract and maybe the Intro some more too

\end{abstract}
 
\vspace{.5cm}
\baselineskip=15pt
\centerline{\bf Key words and phrases}
\medskip 
Pochhammer symbol, generalized hypergeometric series, alternating binomial sum, functional equation

\bigskip
\noindent
{\bf 2010 MSC numbers}
\newline{33C20, 05A10} %--to look up.  including combinatorics wrt binomial sums?

\baselineskip=25pt

\pagebreak
%\medskip
\centerline{\bf Introduction and statement of results}

\medskip
Let $(a)_n=\Gamma(a+n)/\Gamma(a)=(-1)^n\Gamma(1-a)/\Gamma(1-a-n)$ denote the Pochhammer symbol,
$\psi=\Gamma'/\Gamma$ the digamma function, and $_pF_q$ the generalized hypergeometric function \cite{nbs,andrews,bailey,grad,prudvol3}.  The following reductions of particular $_3F_2(1)$ functions have been known for some time (\cite{prudvol3}, p. 537).
\newline{\bf Proposition 1}.
$$_3F_2(a,a+1/3,a+2/3;1/3,2/3;1)=2\cdot 3^{-3a/2-1}\cos {{\pi a} \over 2}, ~~\mbox{Re}~a<0,$$
$$_3F_2(a,a+1/3,a+2/3;2/3,4/3;1)={2 \over {1-3a}}3^{-(3a+1)/2}\cos(3a+1){\pi \over 6}, ~~
\mbox{Re} ~a<1/3,$$
and
$$_3F_2(a,a+1/3,a+2/3;4/3,5/3;1)={4 \over {(1-3a)(2-3a)}}3^{-3a/2}\cos(3a+2){\pi \over 6}, ~~
\mbox{Re} ~a<2/3, ~~a \neq 1/3.$$
However, in the words of one of the surviving authors of \cite{prudvol3}, as to the original
proof, ``it is impossible to find the sources now" \cite{brychkovquote}.
Herein we provide a detailed proof of this Proposition, making use of the properties of closely
related alternating binomial sums.  % how about an application--any particular citation(s) of these
% sums in the literature?
We avoid any use or reliance on the WZ method for hypergeometric summation (e.g., \cite{andrews},
section 3.11).  

There are several known transformations for functions $_3F_2(1)$, as illustrated in
Appendix B.  Therefore the left sides in Proposition 1 may be rewritten in terms of other $_3F_2(1)$
functions with altered parameters.

We define the three binomial sums
$$f_{3j}(a) \equiv \sum_{\ell=0}^\infty (-1)^\ell{{-3a} \choose {3\ell+j}}, ~~j=0,1,2. \eqno(1.1)$$
{\bf Proposition 2} (Common functional equation.)
$$3^6 f_{3j}(a+4)=f_{3j}(a), ~~j=0,1,2.$$

{\bf Proposition 3}.  % To give condition(s) on $a$:
$$_3F_2(a,a+1/3,a+2/3;5/3,7/3;z)$$
$$=\Scale[0.85]{
\frac{8 \left(-\left(1-\sqrt[3]{z}\right)^{3-3 a}-(-1)^{2/3} \left(\sqrt[3]{-1}
   \sqrt[3]{z}+1\right)^{3-3 a}+\sqrt[3]{-1} \left(1-(-1)^{2/3} \sqrt[3]{z}\right)^{3-3 a}+3
   (a-1) \sqrt[3]{z} \left(\left(1-\sqrt[3]{z}\right)^{3-3 a}+\left(\sqrt[3]{-1}
   \sqrt[3]{z}+1\right)^{3-3 a}+\left(1-(-1)^{2/3} \sqrt[3]{z}\right)^{3-3
   a}\right)\right)}{81 (a-1) \left(a-\frac{2}{3}\right) \left(a-\frac{1}{3}\right) (3 a-4)
   z^{4/3}}},$$
giving
$$_3F_2(a,a+1/3,a+2/3;5/3,7/3;-1)=
\frac{8^{2-a} (2-3 a)-8 \sqrt{3} \sin (\pi  a)+8 (6 a-7) \cos (\pi  a)}{9 (a-1) (3 a-4) (3
   a-2) (3 a-1)}
$$
and
$$_3F_2(a,a+1/3,a+2/3;5/3,7/3;1)=
\frac{8\ \cdot 3^{-\frac{3 a}{2}-4} \left(3 \sqrt{3} (6 a-5) \sin \left(\frac{\pi  a}{2}\right)+9
   \cos \left(\frac{\pi  a}{2}\right)\right)}{(a-1) \left(a-\frac{2}{3}\right)
   \left(a-\frac{1}{3}\right) (3 a-4)}.
$$

In fact the proof of Proposition 3 indicates how to obtain a family of values
$_3F_2(a,a+1/3,a+2/3;p/3,q/3;\pm 1)$ from $_3F_2(a,a+1/3,a+2/3;p/3-1,q/3-1;\pm 1)$.
Proposition 3 supplements the following expressions for $_3F_2(-1)$ \cite{prudvol3} (p. 547)
which we restate.  In light of the proofs of Propositions 1 and 3, we forego giving a proof.
% to give:  those 3 $_3F_2(-1)$ relations ...
\newline{\bf Proposition 4}. (\cite{prudvol3})
$$_3F_2(a,a+1/3,a+2/3;1/3,2/3;-1)={2 \over 3}(2^{-3a-1}+\cos a\pi), ~~\mbox{Re}~a<1/3,$$
$$_3F_2(a,a+1/3,a+2/3;2/3,4/3;-1)={2 \over {3(1-3a)}}\left[2^{-3a}+\cos{{(3a+1)} \over 3}\pi
\right], ~~\mbox{Re} ~a<0,$$
and
$$_3F_2(a,a+1/3,a+2/3;4/3,5/3;-1)={4 \over {3(1-3a)(2-3a)}}\left[2^{-3a+1}+\cos{{(3a+2)} \over 3}
\pi \right], ~~\mbox{Re} ~a<-1/3.$$

{\bf Proposition 5}.  (Ordinary differential equation).
The function $u(a)={}_3F_2(a,a+1/3,a+2/3;1/3,2/3;1)$ satisfies the differential equation
$$u''(a)+3(\ln 3) u'(a)+{1 \over 4}(\pi^2+\ln^2 27)u(a)=0.$$

{\bf Corollary}.  We have the identity
$$u'(a)=3\sum_{j=0}^\infty {{(a)_j(a+1/3)_j(a+2/3)_j} \over {(1/3)_j(2/3)_j j!}}[\psi(3(j+a))-
\psi(a)]=-3^{-3a/2-1}\left[\pi \sin {{\pi a} \over 2}+\cos {{\pi a} \over 2} \ln 27\right].$$

Appendix A on Pochhammer symbols and certain of their products complements the proofs and the
rest of the discussion.  % have something on the difficulty to extend to other related $_3F_2(1)$'s?
% --to consider.
These products could be used, for instance, as exercises in a beginning graduate level course on 
special functions or analysis.  

The finite series special cases $f_{30}(-n/3)$ and $f_{31}(-n/3)$ for integer $n \geq 0$ occur in the
online database OEIS in sequences A057681 and A057682 respectively.  Thus for these integer sequences
generating functions are readily available.

%\pagebreak
\centerline{\bf Proof of Propositions}

{\it Proposition 1}.
By using a case from (A.1) of Appendix A for $(1/3)_j(2/3)_j$ together with (A.2),
$${}_3F_2(a,a+1/3,a+2/3;1/3,2/3;1)=\sum_{j=0}^\infty {{(a)_j(a+1/3)_j(a+2/3)_j} \over
{(1/3)_j(2/3)_j j!}}$$
$$=\sum_{j=0}^\infty {{(3a)_{3j}} \over {(3j)!}}=\sum_{j=0}^\infty (-1)^j{{-3a} \choose {3j}}.$$
We may note that
$$\sum_{j=0}^\infty (-1)^j{{-3a} \choose j}=0,$$
unless $a=0$, in which case the sum is $1$, being a special case of the binomial summation
$$\sum_{j=0}^\infty (-1)^jz^j {{-3a} \choose j}=(1-z)^{-3a}.$$
Then 
$$0=\sum_{\ell=0}^\infty (-1)^\ell{{-3a} \choose {3\ell}}-\sum_{\ell=0}^\infty (-1)^\ell{{-3a} 
\choose {3\ell+1}}+\sum_{\ell=0}^\infty (-1)^\ell{{-3a} \choose {3\ell+2}}.$$

We may write this relation as
$$f_{30}(a)-f_{31}(a)+f_{32}(a)=0, \eqno(2.1)$$
by using the definition (1.1).
By using the recurrence of binomial coefficients ${{a+1} \choose {\ell+1}}={a \choose \ell}
+{a \choose {\ell+1}}$, we obtain the following relations:
$$f_{30}(a)=f_{31}(a-1/3)-f_{31}(a),$$
and
$$f_{31}(a)=f_{32}(a-1/3)-f_{32}(a). \eqno(2.2)$$
When combined with (2.1) we then obtain
$$f_{32}(a)=2f_{31}(a)-f_{31}(a-1/3),$$
$$3f_{32}(a)=2f_{32}(a-1/3)-f_{32}(a-2/3)+f_{32}(a-1/3),$$
and
$$f_{32}(a)=f_{30}(a)-f_{30}(a-1/3). \eqno(2.3)$$
Thus a difference of $f_{30}$ gives $f_{32}$, a difference of $f_{32}$ provides $f_{31}$, and
a difference of $f_{31}$ yields $f_{30}$.  All three of these sums have the initial value
$f_{3j}(0)=1$.  

{\it Related $_3F_2(1)$ expressions}.  These other hypergeometric series of Proposition 1 may
be rewritten as:
$$_3F_2(a,a+1/3,a+2/3;2/3,4/3;1)={1 \over {(1-3a)}}\sum_{j=0}^\infty (-1)^j{{-3a+1} \choose {3j+1}}$$
$$={1 \over {(1-3a)}}f_{31}\left(a-{1 \over 3}\right) \eqno(2.4)$$
and
$$_3F_2(a,a+1/3,a+2/3;4/3,5/3;1)={2 \over {(1-3a)(2-3a)}}\sum_{j=0}^\infty (-1)^j{{-3a+2} \choose {3j+2}}$$
$$={2 \over {(1-3a)(2-3a)}}f_{32}\left(a-{2 \over 3}\right). \eqno(2.5)$$

{\it Evaluation of $f_{30}(a)$}.  We will evaluate this binomial sum as a case of the more general
sum 
$$f_{30}(a;z) \equiv \sum_{\ell=0}^\infty {{-3a} \choose {3\ell}}z^\ell.$$
We note that the factorization
$$1+z=(1+z^{1/3})[1-(-1)^{1/3}z^{1/3}][1+(-1)^{2/3}z^{1/3}],$$
implies and is implied by the identity
$$1-(-1)^{1/3}+(-1)^{2/3}=0. \eqno(2.6)$$
Identity (2.6) is very easily verified directly, as
$$-e^{i\pi/3}+e^{2\pi i/3}=-2\cos\left({\pi \over 3}\right)=-1.$$
The identity (2.6) in turn implies the identities for $\ell \equiv 1$ and $2$ (mod $3$)
$$1+(-1)^\ell(-1)^{\ell/3}+(-1)^{2\ell/3}=0. $$ % \eqno(2.7)
We then obtain
$$f_{30}(a=-n;z)={1 \over 3}\sum_{\ell=0}^{3n}{{3n} \choose \ell}\left[z^{\ell/3}+(-1)^\ell (-1)^{\ell/3} z^{\ell/3}+(-1)^{2\ell/3}z^{\ell/3}\right]$$
$$={1 \over 3}\left\{\left(1+z^{1/3}\right)^{3n}+[1-(-1)^{1/3}z^{1/3}]^{3n}+[1+(-1)^{2/3}z^{1/3}]^{3n}
\right\}.$$
Similarly for the generally nonterminating sum, with the same decomposition, we obtain
$$f_{30}(a;z)={1 \over 3} \left\{\left(1+z^{1/3}\right)^{-3a}+[1-(-1)^{1/3}z^{1/3}]^{-3a}
+[1+(-1)^{2/3}z^{1/3}]^{-3a}\right\}.$$
Hence via identity (2.6)
$$f_{30}(a)\equiv f_{30}(a;z=-1)={1 \over 3} \left\{\left[1+(-1)^{1/3}\right]^{-3a} +[1-(-1)^{2/3}]^{-3a}\right\}$$
$$=3^{-3a-1}\left\{[1-(-1)^{2/3}]^{3a}+[1+(-1)^{1/3}]^{3a}\right\}.$$ % by identity (4)
By using the relations $|1-(-1)^{2/3}|=\sqrt{3}=|1+(-1)^{1/3}|$ and
arg$[1-(-1)^{2/3}]=-\pi/6=-\mbox{arg}[1+(-1)^{1/3}]$ we find
$$f_{30}(a)=2\cdot3^{-3a/2-1}\cos\left({{\pi a} \over 2}\right).$$

From (2.3) and then (2.2) we determine that
$$f_{32}(a)=-3^{-3a/2-1}\left[\cos{{\pi a} \over 2}+\sqrt{3}\sin {{\pi a} \over 2}\right]
=-2\cdot 3^{-3a/2-1}\cos(3a-2){\pi \over 6},$$
and
$$f_{31}(a)=3^{-3a/2-1}\left[\cos{{\pi a} \over 2}-\sqrt{3}\sin {{\pi a} \over 2}\right]
=2\cdot 3^{-3a/2-1}\cos(3a+2){\pi \over 6}.$$
From these explicit expressions we may confirm that
$$f_{31}(a)-f_{32}(a)=2\cdot3^{-3a/2-1}\cos\left({{\pi a} \over 2}\right)=f_{30}(a),$$
as in (2.1).  It then follows from (2.4) and (2.5) that
$${}_3F_2(a,a+1/3,a+2/3;2/3,4/3;1)={3^{-(3a+1)/2} \over {1-3a}}\left[\sqrt{3}\cos{{\pi a} \over 2}
-\sin{{\pi a} \over 2}\right]$$
$$={2 \over {1-3a}}3^{-(3a+1)/2}\cos(3a+1){\pi \over 6},$$
and
$${}_3F_2(a,a+1/3,a+2/3;4/3,5/3;1)=2{3^{-3a/2} \over {(1-3a)(2-3a)}}\left[\cos{{\pi a} \over 2}
-\sqrt{3}\sin{{\pi a} \over 2}\right]$$
$$={4 \over {(1-3a)(2-3a)}}3^{-3a/2}\cos (3a+2){\pi \over 6}.$$
% see verification in Mca .nb of 8 May 2015
\qed

{\it Proposition 2}.
We then see that $f_{30}$ satisfies the functional equation
$$3^6 f_{30}(a+4)=f_{30}(a).$$
From (2.3) and then (2.2) it follows that in fact
$$3^6 f_{3j}(a+4)=f_{3j}(a), ~~j=0,1,2.$$
\qed

{\it Remarks}. The initial slopes $f_{3j}'(0)$ may also be found.  For example
$f_{32}'(a)|_{a=0}=(3\ln 3-\sqrt{3}\pi)/6$.

An extension of Proposition 1 would be to consider the following $_3F_2(1)$ function, using
another product of Pochhammer symbols coming from (A.1) together with (A.2).
$${}_3F_2(a,a+1/3,a+2/3;5/3,7/3;1)=\sum_{j=0}^\infty {{(a)_j(a+1/3)_j(a+2/3)_j} \over
{(5/3)_j(7/3)_j j!}}$$
$$=24\sum_{j=0}^\infty {{(3a)_{3j}(j+1)} \over {(3j+4)!}}
=8\sum_{j=0}^\infty (-1)^j{{-3a} \choose {3j}}{1 \over {(3j+4)(3j+2)(3j+1)}}$$
$$=8\sum_{j=0}^\infty (-1)^j{{-3a} \choose {3j}}\left[{1 \over 6}{1 \over {(3j+4)}}-{1 \over 2}
{1 \over {(3j+2)}}+{1 \over 3}{1 \over {(3j+1)}}\right].$$
Here the last sum on the right side may be directly related to $f_{31}(a-1/3)/(1-3a)$ as occurs
in (2.4).  By introducing another summation, the second summation on the right side may be
written in terms of $f_{32}$ as appears in (2.5).  However, using this approach, the first
summation on the right requires several new summations.  Via Proposition 3 we determine this
$_3F_2(1)$ value.

{\it Proposition 3}.  We use the derivative property
$${d \over {dz}}{}_3F_2(a,b,c;d,e;z)={{abc} \over {de}} {}_3F_2(a+1,b+1,c+1;d+1,e+1;z),$$
starting with
$${}_3F_2(a,a+1/3,a+2/3;2/3,4/3;z)$$
$$=\frac{\left(1-\sqrt[3]{z}\right)^{1-3 a}+e^{\frac{2 i \pi }{3}} \left(1-e^{-\frac{2 i \pi
   }{3}} \sqrt[3]{z}\right)^{1-3 a}+e^{-\frac{2 i \pi }{3}} \left(1-e^{\frac{2 i \pi }{3}}
   \sqrt[3]{z}\right)^{1-3 a}}{3 (3 a-1) \sqrt[3]{z}}.
$$
Then shifting $a \to a-1$ we obtain the expression for ${}_3F_2(a,a+1/3,a+2/3;5/3,7/3;z)$
and the two cases for $z=1$ and $z=-1$.  \qed

{\it Proposition 5}.  The differential equation follows from the explicit expression for
$u(a)$ given in Proposition 1. \qed

{\it Corollary}.  Using the explicit expression for $u(a)$ from Proposition 1 and the
relation for $\partial_b (b)_j$ given in Appendix A, the summation identity follows.  \qed

{\it Remarks}.  In regard to Proposition 5, the second order linear differential equation for
$u(a)$ has positive constant coefficients.  As such, it admits a ready physical interpretation
as the equation of a damped harmonic oscillator with damping proportional to 
$(\pi^2+\ln^2 27)/4$ and spring constant proportional to $3 \ln 3$.

Similarly, it may be determined that the function $u_2(a)={}_3F_2(a,a+1/3,a+2/3;2/3,4/3;1)$ satisfies
the differential equation
$$u''(a)+3\left[\ln 3-{2 \over {1-3a}}\right] u'(a)+\left[{1 \over 4}(\pi^2+\ln^2 27)-{{9\ln 3} \over 
{1-3a}}\right]u(a)=0.$$

It is possible to write several integral representations for the binomial sums $f_{3j}(a)$, these
including % would like to get better ones to have insight/another proof of the functional eqn
$$f_{30}(a)={}_3F_2(a,a+1/3,a+2/3;1/3,2/3;1)$$
$$={{\Gamma(2/3)} \over {\Gamma(a+2/3)\Gamma(-a)}}\int_0^1 t^{a-1/3}(1-t)^{-a-1} {}_2F_1\left(a,a+{1 \over 3};
{1 \over 3};t\right)dt$$
$$={1 \over {\Gamma(a+2/3)}}\int_0^\infty e^{-t}t^{a-1/3} {}_2F_2\left(a,a+{1 \over 3};{1 \over 3},{2 \over 3};
t\right)dt.$$
One may also consider to insert integral representations for binomial coefficients into the summands of
$f_{3j}(a)$.  However, it appears difficult to ensure convergence of the resulting expressions with this
approach.

\pagebreak
\centerline{\bf Appendix A:  On certain products of Pochhammer symbols}
% get from my notes
We first recall the relation
$$(a)_{2n}=2^{2n}\left({a \over 2}\right)_n\left({{a+1} \over 2}\right)_n,$$
that follows by applying the duplication formula of the Gamma function,
$$\Gamma(2x)={2^{2x-1} \over \sqrt{\pi}}\Gamma(x)\Gamma\left(x+{1 \over 2}\right).$$

Similarly from the triplication formula of the Gamma function, 
$$\Gamma(3x)={3^{3x-1/2} \over {2\pi}} \Gamma(x)\Gamma\left(x+{1 \over 3}\right)
\Gamma\left(x+{2 \over 3}\right),$$
we obtain
$$(a)_{3n}=3^{3n}\prod_{j=1}^3 \left({{a+j-1} \over 3}\right)_n=3^{3n}\left({a \over 3}\right)_n
\left({{a+1} \over 3}\right)_n\left({{a+2} \over 3}\right)_n. \eqno(A.1)$$
For we have
$$(a)_{3n}={{\Gamma(a+3n)} \over {\Gamma(a)}}
={{2\pi 3^{a+3n-1/2} \Gamma\left({a \over 3}+n\right)\Gamma\left({{a+1} \over 3}+n\right)
\Gamma\left({{a+2} \over 3}+n\right)} \over {2\pi 3^{a-1/2} \Gamma\left({a \over 3}\right)
\Gamma\left({{a+1} \over 3}\right)\Gamma\left({{a+2} \over 3}\right)}}$$
$$=3^{3n} {{\Gamma\left({a \over 3}+n\right)} \over {\Gamma\left({a \over 3}\right)}}
{{\Gamma\left({{a+1} \over 3}+n\right)} \over {\Gamma\left({{a+1} \over 3}\right)}}
{{\Gamma\left({{a+2} \over 3}+n\right)} \over {\Gamma\left({{a+2} \over 3}\right)}}
=3^{3n}\left({a \over 3}\right)_n \left({{a+1} \over 3}\right)_n\left({{a+2} \over 3}\right)_n.$$

Hence we may note that
$$(3a)_{3n}=3^{3n}\left(a\right)_n \left(a+{1 \over 3}\right)_n\left(a+{2 \over 3}\right)_n. \eqno(A.2)$$
Cases following from (A.1) include:
$$\left({1 \over 3}\right)_n\left({2 \over 3}\right)_n={{(3n)!} \over {n!27^n}},$$
$$\left({2 \over 3}\right)_n\left({4 \over 3}\right)_n={{(3n+1)!} \over {n!27^n}},$$
$$\left({4 \over 3}\right)_n\left({5 \over 3}\right)_n={1\over 2}{{(3n+2)!} \over {n!27^n}},$$
$$\left({4 \over 3}\right)_n\left({5 \over 3}\right)_n={1\over 6}{{(3n+3)!} \over {(n+1)!27^n}},$$
$$\left({5 \over 3}\right)_n\left({7 \over 3}\right)_n={1\over {24}}{{(3n+4)!} \over {(n+1)!27^n}},$$
% see others in my notes of 5/5/15, to add:
$$\left({7 \over 3}\right)_n\left({8 \over 3}\right)_n={1\over {5!}}{{(3n+5)!} \over {(n+1)!27^n}},$$
and
$$\left({7 \over 3}\right)_n\left({8 \over 3}\right)_n={2\over {6!}}{{(3n+6)!} \over {(n+2)!27^n}}.$$

We freely make use of the relation ${n \choose k}=(-1)^k (-n)_k/k!$.  Letting $\psi=\Gamma'/
\Gamma$ denote the digamma function, we also have
$${\partial \over {\partial b}}(b)_j=(b)_j[\psi(b+j)-\psi(b)].$$

\pagebreak
\centerline{\bf Appendix B:  Selected transformations of $_3F_2(1)$}

\medskip
The following three transformations \cite{bailey} may be used in rewriting Proposition 1.
$$~_3F_2(a,b,c;d,e;1) = {{\Gamma(e-a-b)\Gamma(e)} \over {\Gamma(e-a)\Gamma(e-b)}}
~_3F_2(a,b,d-c;d,1+a+b-e;1)$$
$$-{{\Gamma(a+b-e)\Gamma(d)\Gamma(e)\Gamma(d+e-a-b-c)} \over {\Gamma(a)\Gamma(b)\Gamma(d-c)
\Gamma(d+e-a-b)}} ~_3F_2(e-a,e-b,d+e-a-b-c;1+e-a-b,d+e-a-b;1).  \eqno(B.1)$$
$$~_3F_2(a,b,c;d,e;1) = {{\Gamma(1+a-d)\Gamma(1+b-d)\Gamma(1+c-d)\Gamma(d)\Gamma(e)} \over
{\Gamma(a)\Gamma(b)\Gamma(c)\Gamma(1+e-d,2-d)}}$$
$$\times ~_3F_2(1+a-d,1+b-d,1+c-d;1+e-d,2-d;1)$$
$$+ {{\Gamma(1+a-d)\Gamma(1+c-d)} \over {\Gamma(1-d)\Gamma(1+a+c-d)}} ~_3F_2(a,c,e-b;1+a+c-d,e;1).
\eqno(B.2)$$
$$~_3F_2(a,b,c;d,e;1) = {{\Gamma(1+a-d)\Gamma(1+b-d)\Gamma(1+c-d)\Gamma(d)\Gamma(e)} \over
{\Gamma(a)\Gamma(b)\Gamma(c)\Gamma(1+e-d,2-d)}}$$
$$\times ~_3F_2(1+a-d,1+b-d,1+c-d;1+e-d,2-d;1)$$
$$+ {{\Gamma(1+a-d)\Gamma(1+b-d)\Gamma(1+c-d)\Gamma(e)} \over {\Gamma(1-d)\Gamma(1+a+b-d)
\Gamma(1+a+c-d)\Gamma(e-a)}}$$
$$\times ~_3F_2(a,1+a-d,1+a+b+c-d-e;1+a+b-d,1+a+c-d;1).  \eqno(B.3)$$

%+other physical realizations of these xforms?
\pagebreak

\end{document}